
\input amstex  
\documentstyle{amsppt}  

\magnification=\magstep1
\hsize=6.5truein
\vsize=9truein

\font \smallrm=cmr10 at 10truept
 at 7truept
 at 10truept 
 at 10truept

\document

\baselineskip=.15truein

\def \N {\Bbb N}

\def \Id {\mathop{\hbox{\rm Id}}\nolimits}
\def \End {\mathop{\hbox{\rm End}}\nolimits}\def \mod {\mathop{\hbox{\rm
mod}}\nolimits}

\def \e {\frak e}
\def \g {\frak g}
\def \kh {k[[h]]}
\def \loongrightarrow {\relbar\joinrel\relbar\joinrel\rightarrow} 

\def \ug {U(\g)}
\def \uhg {U_h(\g)}

\topmatter

\title
   Tressages des groupes de Poisson \`a dual quasitriangulaire  
\endtitle

\author
        Fabio Gavarini${}^\dag$ \, ,  \;  Gilles Halbout${}^\ddag$   
\endauthor 

\leftheadtext{ Fabio Gavarini, Gilles Halbout }  
\rightheadtext{ Tressages des groupes de Poisson \`a dual
quasitriangulaire } 

\affil 
 ${}^\dag$ \! Universit\`a di  \hbox{ Roma ``Tor Vergata'', 
Dipartimento di Matematica }  --  Roma, ITALY   \\ 
  ${}^\ddag$ \! Institut de  \hbox{ Recherche Math\'ematique Avanc\'ee,
ULP--CNRS   --   Strasbourg, } FRANCE  \\ 
\endaffil

\address\hskip-\parindent
  ${}^\dag$ \! Universit\`a degli Studi di Roma ``Tor Vergata''   ---  
Dipartimento di Matematica   \newline 
  Via della Ricerca Scientifica, 1   ---   I-00133 Roma, ITALY   ---  
e-mail: gavarini\@{}mat.uniroma2.it  \newline
     \newline  
  ${}^\ddag$ \! Institut de Recherche Math\'ematique Avanc\'ee   ---  
e-mail:  halbout\@math.u-strasbg.fr   \newline 
  7, rue Ren\'e{} Descartes   ---   67084 STRASBOURG Cedex, FRANCE 
\endaddress

\abstract 
   Let  $ \g $  be a quasitriangular Lie bialgebra over a field  $ k $  of 
characteristic zero, and let  $ \g^* $  be its dual Lie bialgebra.  We
prove that the formal Poisson group  $ F[[\g^*]] $  is a braided Hopf
algebra.  More generally, we prove that if  $ \big( U_h, R \big) $  is
any quasitriangular QUEA, then  $ \Big( {U_h}^{\!\prime},
Ad(R){\big\vert}_{{U_h}^{\!\prime} \otimes {U_h}^{\!\prime}} \Big) $  
--- where  $ {U_h}^{\!\prime} $  is defined by Drinfeld ---   is a
braided QFSHA.  The first result is then just a consequence of the
existence of a quasitriangular quantization  $ (U_h,R) $  of  $ U(\g) $ 
and of the fact that  $ {U_h}^{\!\prime} $  is a quantization of 
$ F[[\g^*]] $.  
\endabstract  

\endtopmatter

  \footnote""{ 1991 {\it Mathematics Subject Classification:} \ 
Primary 17B37, 81R50 }  
  \footnote""{ Le premier auteur a et\'e en partie financ'e par une
bourse du  {\it Consiglio Nazionale delle Ricerche\/}  \, (Italy) }

\vskip20pt

\centerline {\bf  Introduction }

\vskip10pt

   Soit  $ \g $  une big\`ebre de Lie sur un corps  $ k $  de
caract\'eristique z\'ero; notons  $ \g^* $  la big\`ebre de Lie duale
de  $ \g $ et $ F[[\g^*]] $ l'alg\`ebre
des fonctions sur le groupe de Poisson
formel associ\'e \`a  $ \g^* \, $.  Si  $ \g $ est
quasitriangulaire, munie d'une  $ r $--matrice  $ r $,  cela donne \`a 
$ \g $  certaines propri\'et\'es additionnelles.  Une question
se pose alors: quelle nouvelle structure obtient-on sur la big\`ebre duale 
$ \g^* \, $?  Dans ce travail, nous allons montrer que l'alg\`ebre de 
Hopf-Poisson topologique  $ F[[\g^*]] $ est une alg\`ebre tress\'ee (nous
donnerons la d\'efinition plus loin).  Cela avait \'et\'e d\'emontr\'e
pour  $ \, \g = {\frak s}{\frak l}(2,k) $  par Reshetikhin 
(cf.~[Re]), et g\'en\'eralis\'e au cas o\`u  $ \g $  est de
Kac-Moody de type fini (cf.~[G1]) ou de type affine (cf.~[G2]) 
par le premier auteur.
                                             \par 
  Pour d\'emontrer le r\'esultat, nous allons utiliser les quantifications
d'alg\`ebres enveloppantes.  D'apr\`es Etingof-Kazhdan (cf.~[EK]), toute
big\`ebre de Lie admet une quantification  $ \uhg $,  \`a savoir une
alg\`ebre de Hopf (topologique) sur  $ \kh $  dont la sp\'ecialisation
\`a  $ \, h = 0 \, $  est isomorphe \`a  $ \, \ug \, $  comme alg\`ebre de
Hopf co-Poisson; de plus, si  $ \g $  est  quasitriangulaire et  $ r $ 
est sa  $ r $--matrice,  alors il existe une quantification 
$ \uhg $  qui est aussi quasitriangulaire,  en tant 
qu'alg\`ebre de Hopf,  munie d'une  $ R $--matrice  $ \, R_h  \in
\uhg \otimes \uhg \, $  telle que  $ \; R_h \equiv 1 + r \, h
\; \mod\, h^2 \; $ (o\`u l'on  identifie les espaces
vectoriels  $ \, \uhg $ et $ \ug [[h]] \, $). 
                                             \par 
  D'apr\`es Drinfel'd ({\it cf.}~[Dr]), pour toute alg\`ebre
enveloppante universelle quantifi\'ee  $ U $,  on peut d\'efinir une
sous-alg\`ebre de Hopf  $ U' $  telle que, si la limite semi-classique de 
$ U $ est $ \ug$ (avec  $ \g $ une big\`ebre de Lie), alors la limite
semi-classique de  $ U' $  est  $ F[[\g^*]] $.  Dans notre cas, si l'on
consid\`ere  $ {\uhg}' $,  on peut remarquer que la  $ R $--matrice 
n'appartient pas,  {\it a priori},  \`a  $ \, {\uhg}' \otimes {\uhg}'
\, $;  n\'eanmoins, nous prouvons que son action adjointe  $ \;
{\frak R}_h := {\hbox{\rm Ad}}(R_h) : \, \uhg \otimes \uhg
\loongrightarrow \uhg \otimes \uhg \, $,  $ \, x \otimes y \mapsto R_h
\cdot (x \otimes y) \cdot R_h^{\,-1} \, $,  \; stabilise  $ {\uhg}'
\otimes {\uhg}' $, donc induit par sp\'ecialisation un op\'erateur 
$ {\frak R}_0 $  sur  $ \, F[[\g^*]] \otimes F[[\g^*]] \, $.  Enfin,
les propri\'et\'es qui font de  $ R_h $  une  $ R $--matrice 
font de  $ {\frak R}_h $  un op\'erateur de tressage,
donc il en est de m\^eme pour  $ \, {\frak R}_0 \, $:  ainsi, la
paire  $ \big( F[[\g^*]], {\frak R}_0 \big) $  est une alg\`ebre
tress\'ee.   

\vskip15pt

\centerline { REMERCIEMENTS }

\vskip7pt

 Les auteurs tiennent \`a remercier M.~Rosso et
C.~Kassel pour de nombreux entretiens. 

\vskip2,1truecm

\centerline {\bf \S\; 1. \ D\'efinitions et rappels }

\vskip10pt

  {\bf 1.1  Les objects classiques.}  Soit  $ k $  un corps fix\'e de
caract\'eristique zero.  Dans la suite  $ k $  sera le corps de base de
tous les objets   --- alg\`ebres et big\`ebres de Lie, alg\`ebres de
Hopf, etc.~---   que nous introduirons.
                                                  \par  
   Suivant [CP], \S 1.3, nous appellons big\`ebre de Lie une paire 
$ \, (\g, \delta_\g) \, $  o\`u  $ \g $  est une alg\`ebre de Lie et 
$ \, \delta_\g \colon \, \g \rightarrow \g \otimes \g \, $ 
est une application lin\'eaire antisym\'etrique   --- dite cocrochet
de Lie ---   telle que son dual  $ \, \delta_\g^* \colon \, \g^*
\otimes \g^* \rightarrow \g^* \, $  soit un crochet de Lie et que 
$ \delta_g $  elle m\^eme soit un 1-cocycle de  $ \g $  \`a valeurs
dans  $ \g \otimes \g $.  Le dual lin\'eaire  $ \g^* $ de  $ \g $ est alors
\`a son tour une big\`ebre de Lie.  Suivant [CP], \S 2.1.B, nous appelons
big\`ebre de Lie quasitriangulaire une paire  $ \, (\g,r) \, $  telle que 
$ \, r \in \g \otimes \g \, $  soit solution de l'\'equation de
Yang-Baxter classique (CYBE)  $ \, [r_{1{}2},r_{1{}3}] + [r_{1{}2},
r_{2{}3}] + [r_{1{}3},r_{2{}3}] = 0 \, $  dans  $ \, \g \otimes \g
\otimes \g \, $  et  $ \g $  soit une big\`ebre de Lie par rapport
au cocrochet  $ \, \delta = \delta_\g \, $  defini par  $ \,
\delta(x) = [x,r] \, $;  l'\'el\'ement  $ r $  est alors appel\'e 
$ r $--matrice  de  $ \g $.  
                                                  \par  
   Si  $ \g $  est une alg\`ebre de Lie, son alg\`ebre 
enveloppante universelle  $ U(\g) $  est une alg\`ebre de Hopf; 
si de plus  $ \g $  est une big\`ebre de Lie, alors  $ U(\g) $  
est en fait une alg\`ebre de Hopf co-Poisson (cf.~[CP], \S 6.2.A).  
                                                  \par  
   Soit  $ \g $  une alg\`ebre de Lie quelconque: on appelle 
alg\`ebre de fonctions sur le groupe formel associ\'e \`a  $ \g $, 
ou tout simplement groupe formel associ\'e \`a  $ \g $,  l'espace 
$ \, F[[\g]] := {\phantom{\big(} U(\g) \phantom{\big)}}^{\!\!\!*}
\, $  dual lin\'eaire de  $ U(\g) $.  Comme  $ U(\g ) $  est une
alg\`ebre de Hopf, son dual  $ \, F[[\g]] \, $  est une
alg\`ebre de Hopf formelle (suivant [Di], Ch.~1).  Remarquons que si 
$ G $  est un groupe alg\'ebrique connexe d'alg\`ebre de Lie 
$ \g $,  et  $ F[G] $  est l'alg\`ebre de Hopf des fonctions
reguli\`eres sur  $ G $,  et si  $ {\frak m}_e $  est l'id\'eal
maximal dans  $ F[G] $  des fonctions qui s'annulent au point unit\'e 
$ \, e \in G \, $,  alors l'alg\`ebre de Hopf formelle  $ \, F[[\g]] \, $ 
n'est rien d'autre que la compl\'etion  $ {\frak m}_e $--adique  de 
$ F[G] $  (cf.~[On], Ch.~I).  Lorsque, de plus,  $ \g $  est une big\`ebre
de Lie,  $ F[[\g]] $  est en fait une alg\`ebre de Hopf-Poisson (cf.~[CP],
\S 6.2.A) formelle.  

\vskip7pt

  {\bf 1.2  Tressages et quasitriangularit\'e.}  Soit  $ H $  une
alg\`ebre de Hopf dans une cat\'egorie tensorielle  $ ({\Cal A},\otimes) $ 
(cf.~[CP], \S 5):  $ H $  est dite tress\'ee (cf.~[Re], D\'efinition 2)
s'il existe un automorphisme d'alg\`ebre  $ {\frak R} $  de 
$ H \otimes H $, appel\'e op\'erateur de tressage de 
$ H $,  diff\'erent de  la volte $ \; \sigma \colon a \otimes b \mapsto
b \otimes a \; $  et v\'erifiant  
  $$  \displaylines{ 
   {\frak R} \circ \Delta = \Delta^{\hbox{\smallrm op}}  \cr 
   (\Delta \otimes {\Id}) \circ {\frak R} = {\frak R}_{13} \circ
{\frak R}_{23} \circ (\Delta \otimes {\Id}) \; ,  \qquad 
({\Id} \otimes \Delta) \circ {\frak R} = {\frak R}_{13} \circ
{\frak R}_{12} \circ ({\Id} \otimes \Delta)  \cr }  $$ 
o\`u   $ \Delta^{\hbox{\smallrm op}} = \sigma \circ \Delta $  et 
$ {\frak R}_{12} $,  $ {\frak R}_{13} $  et  $ {\frak R}_{23} $  sont
les automorphismes de  $ H \otimes H \otimes H $  d\'efinis par  $ \,
{\frak R}_{12} = {\frak R} \otimes {\Id} \, $,  $ \,
{\frak R}_{23} = {\Id} \otimes {\frak R} \, $,  $ \,
{\frak R}_{13} = (\sigma \otimes {\Id}) \circ ({\Id} \otimes {\frak R})
\circ (\sigma \otimes  {\Id}) \, $.  
                                                  \par  
   Enfin, dans le cas o\`u  $ H $  est, de plus, une alg\`ebre de Hopf
Poisson, nous dirons que cette alg\`ebre  est tress\'ee  en tant
qu'alg\`ebre de Hopf Poisson  si elle est tress\'ee  en tant qu'alg\`ebre
de Hopf ---   par un tressage qui est aussi un automorphisme
d'alg\`ebre de Poisson.  

\vskip3pt

   Si la paire  $ (H, {\frak R}) $  est une alg\`ebre tress\'ee, il
r\'esulte de la d\'efinition que  $ {\frak R} $  v\'erifie l'\'equation
de Yang-Baxter quantique   --- QYBE dans la suite ---   dans 
$ \End(H^{\otimes 3}) $,  \`a savoir  
  $$  {\frak R}_{12} \circ {\frak R}_{13} \circ {\frak R}_{23} 
= {\frak R}_{23} \circ {\frak R}_{13} \circ {\frak R}_{12}  $$
ce qui entra\^\i{}ne que pour tout  $ \, n \in \N \, $  le groupe des
tresses  $ {\Cal B}_n $  agit sur  $ H^{\otimes n} $;  on peut ensuite
alors obtenir des invariants de noeuds, selon la recette donn\'ee en [CP],
\S 15.12.  

\vskip3pt

  Une alg\`ebre de Hopf  $ H $  (dans une categorie tensorielle) est
dite quasitriangulaire (cf. [Dr], [CP]) s'il existe un \'el\'ement
inversible  $ \, R \in H \otimes H \, $,  appel\'e  $ R $--matrice  de 
$ H $,  tel que  
  $$  \displaylines{
   R \cdot \Delta (a) \cdot R^{-1} = {\hbox{\rm Ad}}(R) (\Delta 
(a)) = \Delta^{\hbox{\smallrm op}}(a)  \cr 
   (\Delta \otimes {\Id}) (R) = R_{13} R_{23} \; ,  \qquad 
({\Id} \otimes \Delta) (R) = R_{13} R_{12}  \cr }  $$ 
o\`u  $ \, R_{12}$, $R_{13}$ et $ R_{23}$ sont des
\'el\'ements de  $H^{\otimes 3} $ d\'efinis par  $ \, R_{12}
= R \otimes 1 \, $,  $ \, R_{23} = 1 \otimes R \, $ et  $ \, R_{13} =
(\sigma \otimes {\Id}) (R_{23}) =  ({\Id} \otimes \sigma) (R_{12})
\, $.  Il r\'esulte alors des identit\'es ci-dessus que  $ R $ 
v\'erifie la QYBE dans  $ H^{\otimes 3} $,  i.e.  
  $$  R_{12} R_{13} R_{23} = R_{23} R_{13} R_{12} \; .  $$
Ainsi, les produits tensoriels de  $ H $--modules  sont munis d'une
action du groupe des tresses.  En outre, il est clair que si 
$ (H,R) $  est quasitriangulaire, alors  $ \big( H, {\hbox{\rm 
Ad}}(R) \big) $  est tress\'ee.  

\vskip7pt

  {\bf 1.3  Les objects quantiques.}  Soit  $ {\Cal A} $  la cat\'egorie
dont les objets sont les  $ \kh $--modules  topologiquement libres et
complets au sens  $ h $--adique, et les morphismes sont les
applications  $ \kh $--lin\'eaires continues. Pour tous  $ V $,  $ W $ 
dans  $ {\Cal A} $,  d\'efinissons  $ \, V \otimes W \, $ comme
\'etant  la limite projective des  $ \kh \big/ (h^n) $--modules 
$ \, \big( V / h^n V \big) \otimes_{\kh \big/ (h^n)} \big( W / h^n W \big)
\, $:  cela fait de  $ {\Cal A} $  une cat\'egorie tensorielle (voir [CP]
pour plus de details).  D'apr\`es Drinfel'd (cf.~[Dr]), on appelle
alg\`ebre enveloppante universelle quantifi\'ee   --- QUEA dans la suite
---   toute alg\`ebre de Hopf dans la cat\'egorie  $ {\Cal A} $  dont la
limite semi-classique (i.e. la sp\'ecialisation en  $ \, h = 0 \, $)  est
l'alg\`ebre enveloppante universelle d'une big\`ebre de Lie.  De m\^eme,
on appelle alg\`ebre de Hopf de s\'eries formelles quantiques   --- QFSHA
dans la suite ---   toute alg\`ebre de Hopf dans la cat\'egorie  $ {\Cal
A} $  dont la limite semi-classique est l'alg\`ebre de fonctions sur un 
groupe formel.  

\vskip3pt

   Dans la suite, nous aurons besoin du r\'esultat suivant:  

\vskip7pt

\proclaim{Th\'eor\`eme 1.4}  (cf.~[EK]) Soit  $ \g $  une big\`ebre de
Lie.  Il existe une QUEA  $ \, \uhg \, $  dont la limite semi-classique
est isomorphe \`a  $ \, \ug \, $;  en outre, il existe un isomorphisme
de  $ \kh $--modules  tel que $ \, \uhg \cong \ug [[h]] \, $.   
                                            \hfill\break   
   \indent   De plus, si  $( \g, r) $  est quasitriangulaire,  
alors il existe une QUEA  $ \, \uhg \, $  comme ci-dessus et un
\'el\'ement  $ \, R_h \in U_h(\g) \otimes U_h(\g) \, $ 
tels que  $ \, \big( U_h(\g), R_h \big) \, $  soit une alg\`ebre de Hopf
quasitriangulaire et  $ \, R_h = 1 + r \, h + O\left( h^2 \right) \, $ 
(avec  $ \, O\left( h^2 \right) \in h^2 \cdot H \otimes H \, $).  
$ \square $   
\endproclaim  

\vskip7pt

  {\bf 1.5  Le foncteur de Drinfeld.}  Soit  $ H $  une alg\`ebre de
Hopf sur  $ \kh $.  Pour tout  $ \, n \in \N $,  on d\'efinit 
$ \; \Delta^n \colon \, H \longrightarrow H^{\otimes n} \; $  par  $ \,
\Delta^0 := \epsilon $,  $ \, \Delta^1 := {\Id}_{\scriptscriptstyle
H} $  et  $ \, \Delta^n := \big( \Delta \otimes
{\Id}_{\scriptscriptstyle H}^{\otimes (n-2)} \big) \circ
\Delta^{n-1} \, $  si  $ \, n > 2 $.  Pour tout sous-ensemble
ordonn\'e  $ \, \Sigma = \{i_1, \dots, i_k\} \subseteq \{1, \dots,
n\} \, $  avec  $ \, i_1 < \dots < i_k \, $,  \, on d\'efinit
l'homomorphisme  $ \; j_{\scriptscriptstyle \Sigma} \colon \,
H^{\otimes k} \longrightarrow H^{\otimes n} \; $  par  $ \;
j_{\scriptscriptstyle \Sigma} (a_1 \otimes \cdots \otimes a_k) :=
b_1 \otimes \cdots \otimes b_n \; $  avec  $ \, b_i := 1 \, $  si 
$ \, i \notin \Sigma \, $  et  $ \, b_{i_m} := a_m \, $  pour 
$ \, 1 \leq m \leq k \, $;  on pose alors  $ \; \Delta_\Sigma :=
j_{\scriptscriptstyle \Sigma} \circ \Delta^k \, $.  
On d\'efinit  aussi $ \; \delta_n \colon \, H \longrightarrow H^{\otimes n}
\; $  par  $ \; \delta_n := \sum_{\Sigma \subseteq \{1, \dots, n\}} 
{(-1)}^{n-|\Sigma|} \Delta_\Sigma \, $,  pour tout  $ \, n \in \N_+
\, $,  et plus  g\'eneralement, pour tout  $ \, \Sigma = \{i_1, \dots,
i_k\} \subseteq \{1, \dots, n\} \, $,  avec  $ \, i_1 < \cdots < i_k \, $, 
on d\'efinit  
  $$  \delta_\Sigma := \sum_{\Sigma' \subseteq \Sigma} 
{(-1)}^{|\Sigma|-|\Sigma'|} \, \Delta_{\Sigma'}\; .   \eqno (1.1)  $$ 
En particulier,  $ \, \delta_{\{1, \dots, n\}} = \delta_n \, $. 
Gr\^ace au principe d'inclusion-exclusion, ceci \'equivaut \`a  
  $$  \Delta_\Sigma = \sum_{\Sigma' \subseteq \Sigma }
\delta_{\Sigma'}   \eqno (1.2)  $$  
pour tout  $ \, \Sigma = \{i_1, \dots, i_k\} \subseteq \{1, \dots, n\}
\, $  avec  $ \, i_1 < \cdots < i_k \, $.  Enfin on d\'efinit le
sous-espace   
  $$  H' := \big\{\, a \in H \,\big\vert\, \delta_n(a) \in h^n
H^{\otimes n} \, \big\} \, ,  $$  
de  $ H $  que nous consid\`ererons muni de la topologie
induite.  Nous avons alors le  

\vskip7pt

\proclaim{Th\'eor\`eme 1.6}  (cf.~[Dr], \S 7, ou [G3])  Soit  $ H $  une
alg\`ebre de Hopf dans la cat\'egorie  $ {\Cal A} $.  Alors  $ \, H' \, $ 
est une QFSHA.  De plus, si  $ \, H = U_h(\g) \, $  est une QUEA ayant 
$ U(\g) $  comme limite semi-classique, alors la limite semi-classique de 
$ {\phantom{\big(} U_h(\g) \phantom{\big)}}^{\!\!\prime} $  est 
$ \, F[[\g^*]] \, $.   $ \square $  
\endproclaim  

\vskip2,1truecm

\centerline {\bf \S \; 2. \  Les r\'esultats principaux }

\vskip10pt

   Du point de vue technique, le r\'esultat principal de cet
article concerne le cadre g\'en\'eral des alg\`ebres de Hopf
quasitriangulaires:  

\vskip7pt

\proclaim{Th\'eor\`eme 2.1}  Soit  $ H $  une alg\`ebre de Hopf
quasitriangulaire dans la cat\'egorie  $ {\Cal A} $,  et soit  $ R $  sa 
$ R $--matrice.  Alors, l'automorphisme interieur  $ \, {\hbox{\rm Ad}}
(R) \colon \, H \otimes H \rightarrow H \otimes H \, $  se
restreint en un automorphisme de  $ \, H' \otimes H' $.  La paire  $ \,
\Big( H', \, {\hbox{\rm Ad}}(R) {\big\vert}_{H' \otimes H'} \Big) \, $ 
est donc une alg\`ebre de Hopf tress\'ee dans la cat\'egorie 
$ {\Cal A} $.   $ \square $  
\endproclaim

   La preuve de ce th\'eor\`eme sera donn\'ee dans le paragraphe 3. 
Mais nous pouvons d\'ej\`a en tirer comme cons\'equence le r\'esultat
principal annonc\'e par le titre et dans l'introduction, qui nous donne
une interpr\'etation g\'eometrique de la  $ r $--matrice  classique:  

\vskip7pt

\proclaim{Th\'eor\`eme 2.2}  Soit  $ \g $  une big\`ebre de Lie 
quasitriangulaire.  Alors l'alg\`ebre de Hopf Poisson topologique 
$ \, F[[\g^*]] \, $  est tress\'ee.  En outre, il existe une
quantification de  $ F[[\g^*]] $  qui est une alg\`ebre de Hopf
tress\'ee dont l'op\'erateur de tressage se sp\'ecialise en celui
de  $ F[[\g^*]] $.  
\endproclaim

\demo{Preuve}  Soit  $ r $  la  $ r $--matrice  de  $ \g $.  D'apr\`es
le Th\'eor\`eme 1.4, il existe une QUEA quasitriangulaire  $ \, \big(
\uhg, R_h \big) \, $  dont la limite semi-classique est exactement 
$ \, \big( \ug, \, r \big) \, $  \`a savoir,  $ \, \uhg \big/ h \, \uhg
\cong \ug \, $  et  $ \, (R-1) \big/ h \equiv r  \; \mod \, h \,
\uhg^{\otimes 2} \, $. Par le Th\'eor\`eme 1.6, la limite
semi-classique de  $ \uhg' $  est  $ F[[\g^*]] $.  Soit  $ \, {\frak R}_h
:= {\hbox{\rm Ad}}(R_h) \, $;  le Th\'eor\`eme 2.1 nous assure que 
 $ \, \Big( \!\! {\phantom{\big(} \uhg \phantom{\big)}}^{\!\!\prime},
\, {\frak R}_h{\big\vert}_{{\uhg}' \otimes {\uhg}'} \, \Big) \, $  est
une alg\`ebre de Hopf tress\'ee, donc sa limite semi-classique  $ \,
\bigg( F[[\g^*]], \, \left( {\frak R}_h{\big\vert}_{{\uhg}' \otimes 
{\uhg}'} \right) {\Big\vert}_{h=0} \bigg) \, $  est tress\'ee aussi. 
De plus, comme  $ {\frak R}_h $  est un automor-   
 \eject   
\noindent   phisme d'alg\`ebre et le crochet de Poisson de  $ F[[\g^*]] $ 
est donn\'e par  $ \, \{ a, b \} = {\big( [ \alpha, \beta ] \big/ h
\big)}{\big\vert}_{h=0} \, $  pour tout  $ \, a $,  $ b \in F[[\g^*]] \, $ 
et  $ \, \alpha $,  $ \beta \in \!\!\! {\phantom{\big(} \uhg
\phantom{\big)}}^{\!\!\prime} \, $  tels que  $ \, \alpha{\vert}_{h=0} = a
\, $,  $ \, \beta{\vert}_{h=0} = b \, $,  nous avons que  $ \left( {\frak
R}_h{\big\vert}_{{\uhg}' \otimes {\uhg}'} \right) {\Big\vert}_{h=0} $  est
aussi un automorphisme d'alg\`ebre de Poisson.   $ \square $  
\enddemo  

\vskip7pt

   Le th\'eor\`eme ci-dessus donne donc une interpr\'etation
g\'eom\'etrique de la  $ r $--matrice  d'une big\`ebre de Lie 
quasitriangulaire.  Ce m\^eme r\'esultat avait \'et\'e d\'emontr\'e
pour  $ \, \g = {\frak s}{\frak l}(2,k) $  par Reshetikhin (cf.~[Re]), et
g\'en\'eralis\'e par le premier auteur au cas o\`u  $ \g $  est de
Kac-Moody de type fini (cf.~[G1], o\`u une analyse plus pr\'ecise
est effectu\'ee) ou de type affine (cf.~[G2]).
                                              \par   
   Le Th\'eor\`eme 2.2 a aussi une cons\'equence importante. 
Soient  $ \g $  et  $ \g^* $  comme ci-dessus, soit  $ {\frak R} $  le
tressage de  $ F[[\g^*]] $,  et soit  $ \, \e \, $  l'id\'eal maximal
(unique) de  $ \, F[[\g^* \oplus \g^*]] = F[[\g^*]] \otimes F[[\g^*]]
\, $  (produit tensoriel topologique, selon [Di], Ch.~1).  Puisque
$ {\frak R} $  est un automorphisme d'alg\`ebre,   $ \, {\frak R} (\e)
= \e \, $  et  $ {\frak R} $  induit un automorphisme d'\'espace
vectoriel  $ \, \overline{\frak R} \colon \, \e \big/ \e^2 \rightarrow \e
\big/ \e^2 \, $;  or $ \, \e \big/ \e^2 \cong \g \oplus \g \, $,  donc
puisque  $ {\frak R} $  est aussi un automorphisme d'alg\`ebre de Poisson,
la restriction  $ \overline{\frak R} $  est un automorphisme d'alg\`ebre
de Lie de  $ \, \g \oplus \g = \e \big/ \e^2 \, $;  l'automorphisme 
$ \overline{\frak R} $  h\'erite aussi des autres propi\'et\'es du
tressage  $ {\frak R} $.  Enfin, le dual  $ \, \overline{\frak R}^*
\colon \, \g^* \oplus \g^* \rightarrow \g^* \oplus \g^* \, $  est un
automorphisme de cog\`ebre de Lie de  $ \g^* \oplus \g^* $,  dot\'e lui
aussi de plusieurs autres propriet\'es duales de celles de 
$ \overline{\frak R} $.  En particulier,  $ {\frak R} $, 
$ \overline{\frak R} $  et  $ \overline{\frak R}^* $  sont
solutions de la QYBE. Il existe donc une action du groupe des
tresses  $ {\Cal B}_n $  sur  $ {F[[\g^* \oplus \g^*]]}^{\otimes n} $, 
sur  $ {(\g \oplus \g)}^{\otimes n} $,  et sur  $ {(\g^* \oplus
\g^*)}^{\otimes n} $  $ (n \in \N) $,  dont on peut tirer des
invariants de noeuds (selon [CP], \S 15.12).  
                                              \par   
   De tels automorphismes de  $ \, \g^* \oplus \g^* \, $  et de 
$ \, \g \oplus \g \, $  ont \'et\'es introduits dans [WX], \S 9; leur
construction est li\'ee \`a la  "$ R $--matrice  globale", qui donne
aussi une interpr\'etation g\'eom\'etrique de la  $ r $--matrice 
classique. Il conviendrait alors de comparer nos r\'esultats et
ceux de [WX] et d'\'etudier parall\`element les propriet\'es de
fonctorialit\'e de notre construction: tout cela fera l'objet d'un
article \`a suivre.    

\vskip1,1truecm

\centerline {\bf \S \; 3. \  D\'emonstration du th\'eor\`eme 2.1 }

\vskip10pt

   Dans cette section  $ (H,R) $  sera une alg\`ebre de Hopf
quasitriangulaire comme dans l'\'enonc\'e du Th\'eor\`eme 2.1.  Nous
voulons \'etudier l'action adjointe de  $ R $  sur  $ H \otimes H $,  o\`u
cette derni\`ere est munie de sa structure naturelle d'alg\`ebre de Hopf;
nous noterons par  $ \tilde\Delta $  son coproduit, d\'efini par  $ \; 
\tilde\Delta := \sigma_{2{}3} \circ (\Delta \otimes
{\Id}_{\scriptscriptstyle H} \otimes {\Id}_{\scriptscriptstyle H})
\circ ({\Id}_{\scriptscriptstyle H} \otimes \Delta) \; $   
o\`u  $ \, \sigma_{2{}3} \, $  d\'esigne la volte dans les positions 
$ 2 $  et  $ 3 $.  Nous noterons aussi  $ \, I := 1 \otimes 1 \, $ 
l'unit\'e dans  $ H \otimes H $.  Selon notre d\'efinition du produit
tensoriel en  $ {\Cal A} $,  on a  $ \, {\big( H \otimes H \big)}' = H'
\otimes H' \, $.  Notre but est de montrer que, bien que  $ R $ 
n'appartienne pas forc\'ement \`a  $ {\big( H \otimes H \big)}' $,  son
action adjointe  $ \, a \mapsto R \cdot a \cdot R^{-1} \, $  laisse stable 
$ \, {\big( H \otimes H \big)}' = H' \otimes H' \, $.
                                                  \par
   Posons tout d'abord, pour  $ \, \Sigma = \{i_1, \dots, i_k\} \subseteq
\{1, \dots,n\} \, $,  toujours avec  $ \, i_1 < \cdots < i_k \, $:
  $$  R_\Sigma := R_{2i_1-1,2i_k} R_{2i_1-1,2i_{k-1}} \cdots
R_{2i_1-1,2i_1} R_{2i_2-1,2i_k} \cdots R_{2i_{k-1}-1,2i_k} R_{2i_k-1,2i_1}
\cdots R_{2i_k-1,2i_1}  $$ 
(produit de  $ k^2 $  termes) o\`u  $ \, R_{r,s} :=
j_{\scriptscriptstyle \{r,s\}}(R) \, $,  en d\'efinissant  $ \,
j_{\scriptscriptstyle \{r,s\}} \colon \, H \otimes H \longrightarrow
H^{\otimes 2n} \, $  comme pr\'ec\'edemment.  Nous noterons toujours 
$ |\Sigma| $  pour le cardinal de  $ \Sigma $  (ici  $ \, |\Sigma| =
k \, $).   

\vskip7pt

\proclaim{Lemme 3.1}
   Dans  $ {\big( H \otimes H \big)}^{\otimes n} $,  pour tout  $ \,
\Sigma \subseteq \{1, \dots, n\} $,  on a:  $ \, {\tilde\Delta}_\Sigma(R) =
R_\Sigma \, $.
\endproclaim 

\demo{Preuve}  Sans perdre de g\'en\'eralit\'e, nous d\'emontrerons le 
r\'esultat pour  $ \, \Sigma \! = \! \{1, \dots, n\} $,  i.e.  
  $$  {\tilde\Delta}_{\{1, \dots, n\}}(R) = R_{\{1, \dots, n\}} = R_{1,2n} 
\cdot R_{1,2n-2} \cdots R_{1,2} \cdot R_{3,2n} \cdots R_{2n-3,2} \cdot 
R_{2n-1,2n} \cdots R_{2n-1,2} \; .  $$ 
   \indent   Le r\'esultat est \'evident au rang  $ \, n = 1 \, $. 
Supposons le acquis au rang  $ \, n \, $,  et montrons le au rang  $ \, n
+ 1 \, $;  par  d\'efinition de  $ \tilde\Delta $  et par les
propri\'et\'es de la  $ R $--matrice  on a  
  $$  \eqalign{
   {\tilde\Delta}_{\{1, \dots, n+1 \}}  (\!R)  &  = \left( {\tilde\Delta} 
\otimes {{\Id}_{\scriptscriptstyle H \otimes H}}^{\!\! \otimes {n-1}}
\right) \! \big( {\tilde\Delta}_{\{1, \dots, n\}} (R) \big) = \left(
{\tilde\Delta} \otimes {{\Id}_{\scriptscriptstyle H \otimes H}}^{\!\!
\otimes {n-1}} \right) \! \big( R_{\{1, \dots, n\}} \big)  \cr
   &  = \sigma_{2{}3} (\Delta \otimes {\Id}_{\scriptscriptstyle H}^{\,
\otimes  2n} )\!\!\left( {\Id}_{\scriptscriptstyle H} \otimes \Delta
\otimes  {\Id}_{\scriptscriptstyle H}^{ \otimes 2(n-1)}\! \right) \!\!
(R_{1,2n} \cdots  \hskip-0,16pt  R_{1,2} \cdots  \hskip-0,16pt  R_{3,2}
\cdots  \hskip-0,16pt   R_{2n-1,2})  \cr
   &  = \sigma_{2{}3} \! \left( \Delta \otimes {\Id}_{\scriptscriptstyle
H}^{\, \otimes  2n} \right) (R_{1,2n+1} \cdots R_{1,3} R_{1,2} \cdots
R_{4,3} R_{4,2} \cdots  R_{2n,3} R_{2n,2})  \cr 
   &  = \sigma_{2{}3} (R_{1,2n+2} R_{2,2n+2} \cdots R_{1,4} R_{2,4}
R_{1,3}  R_{2,3} \cdots R_{5,4} R_{5,3} \cdots R_{2n+1,4} R_{2n+1,3}) 
\cr  
   &  = R_{1,2n+2} R_{3,2n+2} \cdots R_{1,4} R_{3,4} \cdot R_{1,2} R_{3,2} 
\cdots R_{5,4} \cdot R_{5,2} \cdots R_{2n+1,4} R_{2n+1,2}  \cr   
   &  = R_{1,2n+2} \cdots R_{1,4} R_{1,2} R_{3,2n+2} \cdots R_{3,4} 
R_{3,2} \cdots R_{5,4} R_{5,2} \cdots R_{2n+1,4} R_{2n+1,2}  \cr
   &  = R_{\{1, \dots, n+1\}} \, ,  \;\;\;  \hbox{q.e.d.}  \quad 
\square  \cr }  $$   
\enddemo  

\vskip7pt

   Dor\'enavant pour tout,  $ \, a $,  $ b \in \N \, $,  nous utiliserons
la notation  $ \, C^a_b \, $  pour d\'esigner l'entier  $ \, {b \choose a}
:= {\, b! \, \over \, a! (b-a)! \,} \in \N \, $.   

\vskip7pt

\proclaim{Lemme 3.2}  Pour tout  $ \, a \in {\big( H \otimes
H \big)}'$,  et pour tout ensemble  $ \, \Sigma \, $  tel que 
$ \, |\Sigma| > i \, $,  on a  
  $$  {\tilde\Delta}_\Sigma(a) = \sum_{\Sigma' \subseteq \Sigma, \;\,
|\Sigma'| \leq i}  \!\! {(-1)}^{i-|\Sigma'|} \, C_{|\Sigma| - 1 - 
|\Sigma'|}^{i-|\Sigma'|} \, {\tilde\Delta}_{\Sigma'}(a) + O \big(
h^{i+1} \big) \, .  $$ 
\endproclaim  

\demo{Preuve}  Il suffit de prouver l'\'enonc\'e pour  $ \, \Sigma = \{1,
\dots, n\} , $  avec  $ \, n > i \, $.  Gr\^ace \`a (1.2), on a  
  $$  \eqalign{ 
   {\tilde\Delta}_{\{1, \dots, n\}}(a)  &  = \sum_{\bar\Sigma \subseteq
\{1, \dots, n\}} \delta_{\bar\Sigma}(a) = \sum_{\bar\Sigma \subseteq \{1,
\dots,  n\},~|\bar\Sigma| \leq i} \! \delta_{\bar\Sigma}(a) + O \big(
h^{i+1} \big)  \cr
   &  = \sum_{\bar\Sigma \subseteq \{1, \dots, n\},~|\bar\Sigma| \leq i} \;
\sum_{\Sigma' \subseteq \bar\Sigma} {(-1)}^{|\bar\Sigma|-|\Sigma'|} \, 
{\tilde\Delta}_{\Sigma'}(a) + O \big( h^{i+1} \big)  \cr
   &  = \sum_{\Sigma' \subseteq \{1, \dots, n\},~|\Sigma'| \leq i} \!\!\!
{\tilde\Delta}_{\Sigma'}(a) \sum_{\Sigma' \subseteq
\bar\Sigma,~|\bar\Sigma| 
\leq i} {(-1)}^{|\bar\Sigma|-|\Sigma'|} + O \big( h^{i+1} \big)  \cr
   &  = \sum_{\Sigma' \subseteq \{1, \dots, n\},~|\Sigma'| \leq i} \!\!\! 
{\tilde\Delta}_{\Sigma'}(a) \, {(-1)}^{i-|\Sigma'|} \,
C^{i-|\Sigma'|}_{n-1-|\Sigma'|} + O \big( h^{i+1} \big) \, ,  \;\;
\hbox{q.e.d.}  \quad  \square  \cr }  $$   
\enddemo

\vskip7pt

   Avant de nous attaquer au r\'esultat principal, il nous faut encore un
petit rappel technique sur les coefficients du bin\^ome: on peut le
prouver facilement en utilisant le d\'eveloppement en s\'erie formelle de 
$ \, {(1-X)}^{-(r+1)} \, $,  \`a savoir  $ \, {(1-X)}^{-(r+1)} =
\sum\limits_{k=0}^{\infty} C_{k+r}^r  X^k \, $.  

\vskip7pt

\proclaim{Lemme 3.3}  Soient  $ \, r $,  $ s $,  $ t \in \N \, $  tels que 
$ \, r < t $.  On a alors les relations suivantes (o\`u l'on pose  $ \,
C_u^v := 0 \, $  si  $ \, v > u \, $):  
                            \hfill\break   
   \centerline{ $ \displaystyle{ (a) \quad  \sum_{d=0}^{t} {(-1)}^{d}
\, C_{d-1}^{r} \, C_{t}^{d}  = -{(-1)}^{r} \; ,  \qquad   (b) \quad 
\sum_{d=0}^{t} {(-1)}^{d} \, C_{d+s}^{r} \, C_{t}^{d}  = 0 \; .  
\;\;\; \square } $ }  
\endproclaim  

\vskip7pt

   Voici enfin le r\'esultat principal de cette section:

\vskip7pt

\proclaim{Proposition 3.4}  Pour tout  $ \, a \in {\big( H \otimes H
\big)}' \, $,  nous avons  $ \; R \, a \, R^{-1} \in {\big( H \otimes H
\big)}' \, $.   
\endproclaim  

\demo{Preuve}  Comme nous devons montrer que  $ \, R \, a \, R^{-1} \, $ 
appartient \`a  $ {\big( H \otimes H \big)}'$,  nous devons
consid\'erer les  termes  $ \delta_n \! \left( R \, a \, R^{-1} \right) $, 
$ \, n \in \N \, $.  Pour cela r\'e\'ecrivons  $ \, \delta_{\{1, \dots,
n\}} \! \left( R \, a \, R^{-1} \right) $  en utilisant le Lemme 3.1 et
le fait que  $  \tilde\Delta $,  et plus g\'en\'eralement 
$ {\tilde\Delta}_{\{i_1, \dots, i_k\}} $  (pour  $ k \leq n $),  est un 
morphisme d'alg\`ebre:  $ \; \delta_{\{1, \dots, n\}} \left( R \, a \,
R^{-1} \right) = \sum\limits_{\Sigma \subseteq \{1, \dots, n\}}
{(-1)}^{n-|\Sigma|} R_\Sigma \, {\tilde\Delta}_\Sigma(a) \,
R^{-1}_\Sigma \; $.  
                                                \par   
   Nous allons d\'emontrer par r\'ecurrence sur  $ i $  que   
  $$  \delta_{\{1, \dots, n\}} \left( R \, a \, R^{-1} \right) = O \big( 
h^{i+1} \big)  \qquad  \text{pour tout} \quad  0 \leq i \leq n-1 \; .  
\eqno (\star)  $$ 
Autrement dit, on verra que tous les termes du d\'eveloppement limit\'e
\`a  l'ordre  $ \, n-1 \, $  sont nuls, donc  $ \, \delta_n \left( R \, a
\, R^{-1} \right) = O(h^n) \, $,  d'o\`u notre \'enonc\'e.
                                                \par   
   Pour  $ \, i = 0 \, $,  on a, pour chaque  $ \Sigma \, $, 
$ \, {\tilde\Delta}_\Sigma(a) = \epsilon(a) I^{\otimes n} + O(h) \, $,
$ \, R_\Sigma = I^{\otimes n} + O(h) \, $,  et aussi 
$ \, R_\Sigma^{-1} = I^{\otimes n} + O(h) \, $,  d'o\`u  $ \;
\delta_{\{1, \dots, n\}} \left( R \, a \, R^{-1} \right) =
\sum\limits_{k=1}^{n} C_n^k {(-1)}^{n-k} \epsilon(a) \, I^{\otimes n}
+ O(h) = O(h) \; $,  donc le r\'esultat  $ (\star) $  est vrai pour 
$ \, i = 0 \, $.  
                                                \par   
   Supposons le r\'esultat  $ (\star) $  acquis pour tout  $ \, i' < i
\, $.  \'Ecrivons les d\'eveloppements  $ h $--adiques  de 
$ R_\Sigma $  et  $ R_\Sigma^{-1} $  sous la forme  $ \; R_\Sigma
= \sum_{\ell=0}^\infty R_\Sigma^{\,(\ell)} \, h^\ell \; $  et  $ \;
R_\Sigma^{-1} = \sum_{m=0}^\infty R_\Sigma^{\,(-m)} \, h^m \, $.  Par
la proposition pr\'ec\'edente, nous avons une approximation de 
$ {\tilde\Delta}_\Sigma(a) $  \`a l'ordre  $ j $:  
  $$  {\tilde\Delta}_\Sigma(a) \, = \sum_{\Sigma' \subseteq 
\Sigma,~|\Sigma'|\leq j} {(-1)}^{j-|\Sigma'|} \, C_{|\Sigma| - 1 -
|\Sigma'|}^{j  - |\Sigma'|} \, {\tilde\Delta}_{\Sigma'}(a) + O \big(
h^{j+1} \big) \; .  $$   
   Nous avons alors l'approximation de  $ \, \delta_{\{1, \dots, n\}}
\left( R \, a \, R^{-1} \right) $  suivante:  
  $$  \displaylines{
   {} \;   \delta_{\{1, \dots, n\}} \left( R \, a \, R^{-1} \right) =
\sum_{\Sigma \subseteq \{1, \dots, n\}} \sum_{\ell + m \leq i} {(-1)}^{n - 
|\Sigma|} \, R_\Sigma^{\,(\ell)} \, {\tilde\Delta}_\Sigma (a) \, 
R_\Sigma^{\,(-m)} \, h^{\ell + m} + O \big( h^{i+1} \big) =   \hfill {\ } 
\cr 
   {} \,   = \sum_{j=0}^i \,\, \sum_{\ell + m = i - j} \Bigg(
\sum_{\Sb  \Sigma \subseteq \{1, \dots, n\}  \\  |\Sigma| > j  \endSb}
\sum_{\Sb  \Sigma' \subseteq \Sigma  \\  |\Sigma'| \leq j  \endSb}  \! 
{(-1)}^{n-|\Sigma|} \, {(-1)}^{j-|\Sigma'|} \, 
C_{|\Sigma|-1-|\Sigma'|}^{j-|\Sigma'|} \, R^{\,(\ell)}_\Sigma \, 
{\tilde\Delta}_{\Sigma'}(a) \, R^{\,(-m)}_\Sigma +   \hfill {\ }  \cr 
   {\ } \hfill   + \sum_{\Sb  \Sigma \subseteq \{1, \dots, n\}  \\ 
|\Sigma| \leq j  \endSb}  {(-1)}^{n-|\Sigma|} \, R^{\,(\ell)}_\Sigma \, 
{\tilde\Delta}_\Sigma(a) \, R^{\,(-m)}_\Sigma \Bigg) \, h^{\ell + m}
+ O \big( h^{i+1} \big) =  \cr   
   {} \; \hfill   = \sum_{j=0}^i \, \sum_{\ell +m+j=i} \, \sum_{\Sb 
\Sigma' \subseteq \{1, \dots, n\}  \\  |\Sigma'| \leq  j  \endSb}  \!\! 
\Bigg( \sum_{\Sb  \Sigma \subseteq \{1, \dots, n\}  \\  \Sigma' \subseteq
\Sigma, {\ }  |\Sigma| > j  \endSb}  \hskip-14pt  {(-1)}^{n-|\Sigma|} \,
{(-1)}^{j-|\Sigma'|} \, C_{|\Sigma|-1-|\Sigma'|}^{j-|\Sigma'|} \,
R^{\,(\ell)}_\Sigma {\tilde\Delta}_{\Sigma'}(a) \, R^{\,(-m)}_\Sigma +  
{\ }  \cr 
   {\ } \hfill   + {(-1)}^{n-|\Sigma'|} \, R^{\,(\ell)}_{\Sigma'} \, 
{\tilde\Delta}_{\Sigma'}(a) \, R^{\,(-m)}_{\Sigma'} \Bigg) \, h^{\ell + m}
+ O \big( h^{i+1} \big) \; .  \cr }  $$ 
   \indent   Nous noterons (E) la derni\`ere expression entre
parenth\`ese, et  nous montrerons que cette expression est nulle, d'o\`u 
$ \, \delta_n \left( R \, a \, R^{-1} \right) = O \big( h^{i+1} \big)
\, $.   
                                                \par   
   Regardons d'abord les termes correspondant \`a  $ \, \ell + m = 0
\, $,  c'est-\`a-dire  $ \, j = i \, $.  On retrouve alors  $ \delta_{\{1,
\dots, n\}}(a) $,  qui est dans  $ O \big( h^{i+1} \big) $  par
hypoth\`ese.  Dans la suite du calcul nous supposerons d\'esormais 
$ \, \ell + m > 0 \, $.  
 \eject
   Regardons maintenant comment les termes  $ R_\Sigma^{\,(\ell)} $  et 
$ R_\Sigma^{\,(-m)} $  agissent sur  $ {{\big( H \otimes H \big)}'}^{\,
\otimes n} $  (respectivement \`a gauche et \`a droite) pour  $ \, \ell +
m \, $  fix\'e (et positif), disons  $ \, \ell + m = S \, $.  En faisant
le d\'eveloppement limit\'e de chaque  $ R_{i,j} $  qui appara\^\i{}t
dans  $ R_\Sigma \, $,  on voit que  $ R_\Sigma^{\,(\ell)} $  et 
$ R_\Sigma^{\,(-m)} $  sont sommes de produits d'au plus  $ \ell $ 
et  $ m $  termes respectivement, chacun agissant sur au plus deux
facteurs tensoriels de  $ {{\big( H \otimes H \big)}'}^{\, \otimes n} $. 
Nous allons r\'e\'ecrire  $ \, \sum\limits_{\ell + m = S}
R^{\,(\ell)}_\Sigma \, {\tilde\Delta}_{\Sigma'}(a) \,
R^{\,(-m)}_\Sigma \, $  en regroupant les termes de la somme
qui agissent sur les m\^emes facteurs de  $ {{\big( H \otimes H 
\big)}'}^{\, \otimes n} $,  facteurs dont nous identifierons les
positions par  $ \, \Sigma'' $.   
                                             \par   
   Si  $ i $  appartient \`a  $ \Sigma'' $,  dans l'identification 
$ \, {(H \otimes H)}^{\otimes n} = H^{\otimes 2n} \, $  (telle
qu'on l'a choisie pour d\'efinir  $ R_\Sigma \, $)  l'indice  $ i $ 
correspond \`a la paire  $ (2i-1,2i) \, $;  mais alors  $ R_\Sigma $  et 
$ R_\Sigma^{\,-1} $,  et  donc aussi chaque  $ R_\Sigma^{\,(\ell)} $  et
chaque  $ R_\Sigma^{\,(-m)} \, $,  n'agissent de mani\`ere non triviale sur
le  $ i $--\`eme  facteur de  $ {\tilde\Delta}_{\Sigma'}(a) $  que si, dans
l'\'ecriture explicite de  $ R_\Sigma $, un terme non trivial apparait aux
places  $ \, 2i-1 \, $  ou  $ 2i $,  donc seulement si  $ \, i \in \Sigma
\, $:  ainsi  $ \, \Sigma'' \subseteq \Sigma \, $.  Nous posons alors  
  $$  \sum_{\ell + m = S} R^{\,(\ell)}_\Sigma \,
{\tilde\Delta}_{\Sigma'}(a) \, R^{\,(-m)}_\Sigma =
\sum_{\Sigma'' \subseteq \Sigma} A^{(S)}_{\Sigma', \Sigma,
\Sigma''}(a) \, .  $$  
   \indent   Maintenant consid\'erons  $ \, \bar\Sigma \supseteq \Sigma
\, $.  D'apr\`es la d\'efinition on a  $ \, R_{\bar\Sigma} = R_\Sigma +
{\Cal A} \, $,  o\`u  $ {\Cal A} $  est une somme de termes qui
contiennent des facteurs  $ \, R_{2i-1,2j}^{\,(s)} \, $  avec  $ \, \{i,j\}
\not\subseteq \Sigma \, $:  pour ce voir, il suffit de d\'evelopper chaque
facteur  $ R_{a,b} $  dans  $ R_{\bar\Sigma} $  comme  $ \, R_{a,b} =
1^{\otimes 2n} + O(h) \, $.  De  m\^eme, on a aussi  $ \,
R_{\bar\Sigma}^{\,(\ell)} = R_\Sigma^{\,(\ell)} + {\Cal A}' \, $,  et
pareillement  $ \, R_{\bar\Sigma}^{\,(-m)} = R_\Sigma^{\,(-m)} + {\Cal A}''
\, $. Cela implique que  $ \,  A^{(S)}_{\Sigma'', \bar\Sigma, \Sigma'}(a) =
A^{(S)}_{\Sigma'', \Sigma, \Sigma'}(a) \, $,  et donc les 
$ A^{(S)}_{\Sigma'', \Sigma, \Sigma'}(a) $  ne d\'ependent pas de 
$ \Sigma \, $;  on \'ecrit alors  
  $$  \sum_{\ell + m = S} R^{\,(\ell)}_\Sigma \,
{\tilde\Delta}_{\Sigma'}(a) \,  R^{\,(-m)}_\Sigma = \sum_{\Sigma''
\subseteq \Sigma}  A^{(S)}_{\Sigma',\Sigma''}(a) \; .  $$  
   \indent   Nous allons ensuite r\'e\'ecrire  $ (E) $  \`a l'aide des 
$ A^{(S)}_{\Sigma', \Sigma''}(a) $.  Par commodit\'e dans la suite des 
calculs, on notera  $ \, \delta_{\Sigma'' \subseteq \Sigma'} \, $  la
fonction  qui vaut  $ 1 $  si  $ \, \Sigma'' \subseteq \Sigma' \, $ 
et  $ 0 $  sinon.  Nous obtenons alors une nouvelle expression pour 
$ \, \delta_{\{1, \dots,  n\}} \left( R \, a \, R^{-1} \right) \, $, 
\`a savoir  
  $$  \displaylines{
   {} \;   \delta_{\{1, \dots, n\}} \! \left( R \, a \, R^{-1} \right)
= \sum_{j=0}^{i-1} \, \sum_{\Sb  \Sigma' \subseteq \{1, \dots, n\}  \\ 
   |\Sigma'| \leq  j  \endSb} \Bigg( \sum_{\Sb  \Sigma \subseteq \{1,
\dots, n\}  \\  \Sigma' \subseteq \Sigma, {\ } |\Sigma| > j  \endSb} 
\!\!\!\! {(-1)}^{n-|\Sigma|} \, {(-1)}^{j-|\Sigma'|} \, 
C_{|\Sigma|-1-|\Sigma'|}^{j-|\Sigma'|} \times  \cr   
   {\ } \hfill   \times \sum_{\Sigma'' \subseteq \Sigma}
\! A^{(i-j)}_{\Sigma', \Sigma''}(a) + {(-1)}^{n-|\Sigma'|}
\sum_{\Sigma'' \subseteq \Sigma'} A^{(i-j)}_{\Sigma', \Sigma''}(a)
\Bigg) \, h^{i-j} + O \big( h^{i+1} \big) = {\ }  \cr   
   {\;} = \, \sum_{j=0}^{i-1} \, \sum_{\Sb  \Sigma' \subseteq 
\{1, \dots, n\}  \\  |\Sigma'| \leq  j  \endSb}  h^{i-j}  \sum_{\Sigma'' 
\subseteq \{1, \dots, n\}} \! A^{(i-j)}_{\Sigma', \Sigma''}(a) \times
\hfill {\ }  \cr 
   {\ } \hfill   \times \Bigg( \sum_{\Sb  \Sigma \subseteq \{1, \dots,
n\}  \\  
\Sigma' \subseteq \Sigma, \; \Sigma'' \subseteq \Sigma, \; |\Sigma| > j
\endSb}  \!\!\!\!\! {(-1)}^{n-|\Sigma|} \, {(-1)}^{j-|\Sigma'|} \,
C_{|\Sigma| -  1 - |\Sigma'|}^{j - |\Sigma'|} + {(-1)}^{n-|\Sigma'|} \,
\delta_{\Sigma'' \subseteq \Sigma'} \Bigg)  + O \big( h^{i+1} \big) \; . 
{\ }  \cr }  $$ 
   \indent   Notons  $ {\big( E' \big)}_{\Sigma',\Sigma''} $  
la nouvelle expression entre parenth\`ese; autrement dit, pour 
$ \Sigma' $  et  $ \Sigma'' $  fix\'ees, avec  $ \, \big| \Sigma' \big|
\leq j \, $,  on pose   
  $$  {\big( E' \big)}_{\Sigma',\Sigma''} := \sum_{\Sb  \Sigma \subseteq
\{1, \dots, n\}  \\  \Sigma' \subseteq \Sigma, \; \Sigma'' \subseteq
\Sigma, \; |\Sigma| > j  \endSb} \!\!\!\!\! {(-1)}^{n-|\Sigma|} \,
{(-1)}^{j-|\Sigma'|} \, C_{|\Sigma| - 1 -  |\Sigma'|}^{j - |\Sigma'|} +
{(-1)}^{n-|\Sigma'|} \, \delta_{\Sigma'' \subseteq \Sigma'}  $$  
(au passage, remarquons que celle-ci est une expression purement
combinatoire); nous allons montrer que cette expression est nulle lorsque 
$ \Sigma' $  et  $ \Sigma'' $  sont telles que  $ \, \big| \Sigma' \cup
\Sigma'' \big| \leq j - i + \big| \Sigma' \big| \, $ et $\,\big| \Sigma'
\big| \leq j \,$.  En vertu du  lemme suivant, ceci suffira pour
prouver la proposition.   

\vskip7pt

\proclaim{Lemme 3.5}
                                          \hfill\break
   \indent   {\ } \; {\it (a)} \;  On a  $ \, j < i \, $  et 
$ \, i \leq n-1 \, $,  donc  $ \, j \leq n-2 \, $.
                                          \hfill\break
   \indent   {\ } \; {\it (b)} \;  Pour tout  $ \, S > 0 \, $,  dans 
l'expression  $ \; \sum\limits_{\ell + m = S} R^{\,(\ell)}_\Sigma \, 
{\tilde\Delta}_{\Sigma'}(a) \, R^{\,(-m)}_\Sigma = \sum\limits_{\Sigma'' 
\subseteq \Sigma} A^{(S)}_{\Sigma',\Sigma''}(a) \, $  \; on a  
$ \; A^{(S)}_{\Sigma', \Sigma''}(a) = 0 $  \ pour tout  $ \Sigma' $, 
$ \Sigma'' $  tels que  $ \, \big| \Sigma' \cup \Sigma'' \big| > S + 
\big| \Sigma' \big| \, $.
\endproclaim

\demo{Preuve}  La premi\`ere assertion est \'evidente; pour montrer la
deuxi\`eme nous \'etudions l'action adjointe de  $ R_\Sigma $  sur 
$ {\big( H \otimes H \big)}^{\, \otimes n} $.
                                          \par
   Premi\`erement, sur  $ k \cdot I^{\otimes n} $  l'action de ces
\'el\'ements  donne un terme nul car on retrouve le terme \`a l'ordre  $ S
$  du  d\'eveloppement  $ h $--adique  de  $ \, R_\Sigma \cdot
R_\Sigma^{-1} = 1 \, $   (pour  $ \, S > 0 \, $).
                                          \par
   Deuxi\`emement, consid\'erons  $ \, \Sigma \subseteq \{1, \dots, n\}
\, $,   et \'etudions l'action sur  $ \, {\big( H \otimes H)}_{\Sigma'} := 
j_{\scriptscriptstyle \Sigma'} \left( {\big( H \otimes H \big)}^{\otimes 
|\Sigma|} \right) \, (\, \subseteq {\big( H \otimes H \big)}^{\otimes n}
\,) \, $.  On sait que  $ R_\Sigma $  est un produit de  $ {|\Sigma|}^2 $ 
termes du  type  $ \, R_{a,b} \, $,  avec  $ \, a, b \in \big\{\, 2i-1, 2j
\,\big\vert\, i,  j \in \Sigma \,\big\} \, $;  analysons donc ce qui se
passe lorsqu'on fait le  produit  $ \, P := R_\Sigma \cdot x \cdot
R_\Sigma^{\,-1} \, $  si  $ \, x \in  {\big( H \otimes H)}_\Sigma \, $.
                                         \par
   Consid\'erons le facteur  $ R_{a,b} $  qui appara\^\i{}t le plus \`a
droite: si  $ \, a, b \not\in \big\{\, 2j-1, 2j \,\big\vert\, j \in
\Sigma' \,\big\} \, $,  alors en calculant  $ P $  on trouve  $ \, P :=
R_\Sigma \, x \,  R_\Sigma^{\,-1} = R_\star \, R_{a,b} \, x \,
R_{a,b}^{\,-1} \, R_\star^{\,-1} =  R_\star \, x \, R_\star^{\,-1} \, $ 
(o\`u  $ \, R_\star := R_\Sigma \,  R_{a,b}^{\,-1} \, $).  De m\^eme, en
avan\c{c}ant de droite \`a gauche le long  de  $ R_\Sigma $  on peut
\'ecarter tous les facteurs  $ R_{c,d} $  de ce type, \`a savoir tels que 
$ \, c, d \not\in \big\{\, 2j-1, 2j \,\big\vert\, j \in \Sigma' \,\big\}
\, $.  Ainsi le premier facteur dont l'action adjointe est non  triviale
sera n\'ecessairement du type  $ R_{\bar{a},\bar{b}} $  avec l'un des 
deux indices appartenant \`a  $ \, \big\{\, 2j-1, 2j \,\big\vert\, j
\in \Sigma' \,\big\} \, $,  soit par exemple  $ \bar{a} $.  Notons que le
nouvel indice  $ \, \bar{a} \; (\, \in \{1, 2, \dots, 2n-1, 2n\} \,)
\, $,  qui  agit sur un facteur tensoriel dans  $ H^{\otimes 2n} $, 
correspond \`a un nouvel indice  $ \, j_{\bar{a}} \; (\, \in \{1, \dots,
n\} \,) \, $,  agissant  sur un facteur tensoriel de  $ {\big( H \otimes
H \big)}^{\otimes n} \, $.  Ainsi pour les facteurs successifs   ---
i.e.~\`a gauche de  $ R_{\bar{a},\bar{b}} $  ---   il faut r\'ep\'eter la
m\^eme analyse, mais avec l'ensemble  $ \, \big\{\, 2j-1, 2j \,\big\vert\,
j \in \Sigma' \cup \{ j_{\bar{a}} \} \,\big\} \, $  \`a la place de  $ \,
\big\{\, 2j-1, 2j \,\big\vert\, j \in \Sigma' \,\big\} \, $;  donc, comme 
$ R_{\bar{a},\bar{b}} $  pouvait agir de mani\`ere non triviale sur au
plus  $ \big| \Sigma' \big| $   facteurs de  $ {\big( H \otimes H
\big)}^{\otimes n} $,  de m\^eme le facteur  plus proche \`a sa gauche ne
peut agir de mani\`ere non triviale que sur au plus  $ \, \big| \Sigma'
\big| + 1 \, $  facteurs.  La conclusion est que  l'action adjointe
de  $ R_\Sigma $  est non triviale sur au plus  $ \, \big| \Sigma'
\big| + \big| \Sigma \big| \, $  facteurs de 
$ {\big( H \otimes H \big)}^{\otimes n} $.
                                          \par
   Maintenant, consid\'erons les diff\'erents termes 
$ R_\Sigma^{\,(\ell)} $   et  $ R_\Sigma^{\,(-m)} \,$, 
avec  $ \, \ell + m = S \, $,  et \'etudions les  produits 
$ \, R_\Sigma^{\,(\ell)} \cdot x \cdot R_\Sigma^{\,(-m)} \, $, 
avec  $ \, x \in {\big( H \otimes H)}_\Sigma \, $.  On sait d\'ej\`a que 
$ R_\Sigma^{\,(\ell)} $  et  $ R_\Sigma^{\,(-m)} \, $  sont sommes de
produits,  not\'es  $ P_+ $  et  $ P_- \, $,  d'au plus  $ \ell $  et 
$ m $  termes  respectivement, du type  $ \, R_{i,j}^{\,(\pm k)} \, $; 
les termes  $ A_{\Sigma',\Sigma''}^{(S)}(a) $  ne sont alors que des sommes
de termes du type  $ \, P_+ \, {\tilde\Delta}_{\Sigma'}(a) \, P_- \, $,
o\`u de plus les ``indices'' intervenant dans  $ P_+ $  et  $ P_- \, $ 
sont dans  $ \Sigma'' $.  Or,  comme chaque  $ P_+ $  et  chaque  $ P_- $ 
est un produit d'au plus  $ \ell $   et  $ m $  facteurs  $ \,
R_{i,j}^{\,(\pm k)} \, $,  on peut raffiner l'argument  pr\'ec\'edent. 
Consid\'erons seulement le terme \`a l'ordre  $ S $  du d\'eveloppement 
$ h $--adique  de  $ \, P := R_\Sigma \, x \,  R_\Sigma^{\,-1} = R_\star
\, R_{a,b} \, x \, R_{a,b}^{\,-1} \, R_\star^{\,-1} =  R_\star \, x \,
R_\star^{\,-1} \, $:  lorsque il y a des facteurs du type 
$ R_{a,b}^{\,(k)} $  ou  $ R_{a,b}^{\,(t)} \, $,  pour  $ a $,  $ b $  
fix\'es   --- n'appartenant pas \`a  $ \big\{\, 2j-1, 2j \,\big\vert\,
j \in \Sigma' \,\big\} \, $  ---   qui appara\^\i{}ssent dans 
$ R_\Sigma^{\,(\ell)} $   ou  $ R_\Sigma^{\,(-m)} \, $,  pour certains 
$ \ell $  ou  $ m $,  la contribution totale de tous ces termes dans la
somme  $ \, \sum\limits_{\ell + m = S} R_\Sigma^{\,(\ell)} \, x \,
R_\Sigma^{\,(-m)} \, $   sera nulle (cela vient du fait que  $ R_\star
\, R_{a,b} \, x \, R_{a,b}^{\,-1} \, R_\star^{\,-1} = R_\star \, x \,
R_\star^{\,-1} \, $).  De plus, comme maintenant on ne consid\`ere que 
$ S $  facteurs au total, on conclut que  $ \; A_{\Sigma',
\Sigma''}^{(S)}(a) = 0 \; $  si  $ \; \big| \Sigma' \cup 
\Sigma'' \big| > S + \big| \Sigma' \big| \, $.    
                                      $ \square $\break      
\enddemo

\vskip5pt

   Calculons maintenant  $ {\big( E' \big)}_{\Sigma',\Sigma''} $.  
Gr\^ace \`a la remarque pr\'ec\'edente, nous pouvons nous limiter aux
paires  $ \big( \Sigma', \Sigma'' \big) $  telles que  $ \; \big| \Sigma'
\cup \Sigma'' \big| \leq i - j + m + \big| \Sigma' \big| \leq i - j + j =
i \leq n-1 \; $.  On pourra alors toujours trouver au moins deux  $ \,
\Sigma \subseteq \{1, \dots, n\} \, $  tels que  $ \, |\Sigma| > j \, $ 
et  $ \, \Sigma' \cup \Sigma'' \subseteq \Sigma \, $,  ce qui nous assure
qu'il y'aura toujours au  moins deux termes dans le comptage qui va suivre
(condition qui assurera la  nullit\'e de l'expression  $ {\big( E'
\big)}_{\Sigma',\Sigma''} \, $).  Nous  allons distinguer trois cas:

\vskip3pt

 (I) \quad  Si  $ \, \Sigma'' \subseteq \Sigma' \, $,  alors l'expression  
$ {\big( E' \big)}_{\Sigma',\Sigma''} $  devient  
  $$  {\big( E':1 \big)}_{\Sigma',\Sigma''}
 \, = \sum_{\Sb  \Sigma \subseteq \{1, \dots, n\}  \\  
\Sigma' \subseteq \Sigma, \; |\Sigma| > j  \endSb}  \!\!
{(-1)}^{n-|\Sigma|} \,  {(-1)}^{j-|\Sigma'|} \,
C_{|\Sigma|-1-|\Sigma'|}^{j-|\Sigma'|} +  {(-1)}^{n-|\Sigma'|} \; .  $$  
En regroupant les  $ \Sigma $  qui ont le m\^eme cardinal  $ d $,  un
simple  comptage nous donne  
  $$  { \big( E':1 \big)}_{\Sigma',\Sigma''} = \sum_{d=j+1}^n {(-1)}^{n-d}
\,  {(-1)}^{j-|\Sigma'|} 
\, C_{d-1-|\Sigma'|}^{j-|\Sigma'|} \, C_{n-|\Sigma'|}^{d-|\Sigma'|} + 
{(-1)}^{n-|\Sigma'|} \; .  $$  
   \indent   Or, cette derni\`ere expression est nulle d'apr\`es le
Lemme 3.3, car elle correspond \`a une somme du type  $ \;
\sum\limits_{k=r+1}^t {(-1)}^{t+r-k} \, C_{k-1}^r \, C_t^k + {(-1)}^t =
\sum\limits_{k=0}^t {(-1)}^{t+r-k} \, C_{k-1}^r \, C_t^k + {(-1)}^t \; $ 
(o\`u  $ \, C_u^v := 0 \, $  si  $ \, v > u \, $)  avec  $ \, r $, 
$ t \in \N_+ \, $  et  $ \, r < t \, $:  dans notre cas on a pos\'e  $ t =
n - \big| \Sigma' \big| \, $,  $ \, r = j - \big| \Sigma' \big| \, $  et 
$ \, k = d - \big| \Sigma' \big| \, $;  on v\'erifie que l'on a  $ \, j
- \big| \Sigma' \big| < n - \big| \Sigma' \big| \, $  parce que  $ \, j <
n \, $.

\vskip3pt

 (II) \quad  Si  $ \, \Sigma'' \not\subseteq \Sigma' \, $  et 
$ \, \big| \Sigma' \cup \Sigma'' \big| > j \, $,  alors l'expression  $
{\big(  E' \big) }_{\Sigma',\Sigma''}$  devient   
  $$  {\big( E':2 \big)}_{\Sigma',\Sigma''}
\, = \sum_{\Sb  \Sigma \subseteq \{1, \dots, n\}  \\ 
\Sigma' \cup \Sigma'' \subseteq \Sigma  \endSb}  \!\! {(-1)}^{n-|\Sigma|}
\,  {(-1)}^{j-|\Sigma'|} \, C_{|\Sigma|-1-|\Sigma'|}^{j-|\Sigma'|} \; .  $$
En regroupant les  $ \Sigma $  qui ont le m\^eme cardinal  $ d $,  un
simple  comptage nous donne  
  $$  {\big( E':2 \big) }_{\Sigma',\Sigma''}\, = \sum_{d = |\Sigma' \cup 
\Sigma''|}^n
{(-1)}^{n-d} \, {(-1)}^{j-|\Sigma'|} \, C_{d-1-|\Sigma'|}^{j-|\Sigma'|} \,
C_{n - |\Sigma' \cup \Sigma''|}^{d - |\Sigma' \cup \Sigma''|} \; .  $$
   \indent   \`A nouveau, cette derni\`ere expression est nulle gr\^ace
au Lemme 3.3, car elle correspond \`a une somme du type  $ \;
\sum\limits_{k=0}^t {(-1)}^{t+r-k} \, C_{k+s}^r \, C_t^k \; $ 
avec  $ \, r $,  $ t $,  $ s \in \N_+ \, $  et 
                      $ \, r < t \, $:  dans notre cas\break     
 \eject  
\noindent   on a  pos\'e  $ \, t = n - \big|
\Sigma'\cup \Sigma'' \big| \, $,  $ \, r = j - \big| \Sigma' \big| \, $, 
$ \, s = \big| \Sigma' \cup \Sigma'' \big| - \big| \Sigma' \big| - 1 \, $ 
et  $ \, k = d - \big| \Sigma' \cup \Sigma'' \big| \, $;   on v\'erifie que
l'on a  $ \, j - \big| \Sigma' \big| < n - \big| \Sigma' \big| \, $  car 
$ \, j < n \, $  et  $ \big| \Sigma' \cup \Sigma'' \big| - \big| \Sigma'
\big| - 1 \geq 0 \, $  car  $ \, \Sigma'' \not\subseteq \Sigma' \, $.

\vskip3pt

 (III) \quad  Si  $ \, \Sigma'' \not\subseteq \Sigma' \, $  et 
$ \big| \Sigma' \cup \Sigma'' \big| \leq j \, $,  alors l'expression 
$ {\big( E' \big)}_{\Sigma',\Sigma''} $  devient  
  $$  {\big( E':3 \big)}_{\Sigma',\Sigma''}
 \, = \sum_{\Sb  \Sigma \subseteq \{1, \dots, n\} \\  
\Sigma' \cup \Sigma'' \subseteq \Sigma, \; |\Sigma|>j  \endSb}  
{(-1)}^{n-|\Sigma|} \, {(-1)}^{j-|\Sigma'|} \, 
C_{|\Sigma|-1-|\Sigma'|}^{j-|\Sigma'|} \; .  $$ 
Si l'on regroupe les  $ \Sigma $  qui ont le m\^eme cardinal  $ d $,  un
simple  comptage nous donne  
  $$ { \big( E':3 \big)}_{\Sigma',\Sigma''}
 = \sum_{d=j+1}^n {(-1)}^{n-d} \, {(-1)}^{j-|\Sigma'|} 
\, C_{d-1-|\Sigma'|}^{j-|\Sigma'|} \, C_{n - |\Sigma' \cup \Sigma''|}^{d - 
|\Sigma' \cup \Sigma''|} \; .  $$ 
   \indent   Mais encore la derni\`ere expression est nulle d'apr\`es le
Lemme 3.3, car elle correspond \`a une somme du type  $ \; \sum\limits_{k
= j + 1 - |\Sigma' \cup \Sigma''|}^t \!\!  {(-1)}^{t+r-k} \, C_{k+s}^r \,
C_t^k = \sum\limits_{k=0}^t {(-1)}^{t+r-k} \,  C_{k+s}^r \, C_t^k \; $ 
(o\`u  $ \, C_u^v := 0 \, $  si  $ \, v > u \, $)  avec  $ \, r $,  $ t $, 
$ s \in \N_+ \, $  et  $ \, r < t \, $:  ici on a encore pos\'e  $ \, t = n
- \big| \Sigma' \cup \Sigma'' \big| \, $,  $ \, r = j - \big| \Sigma'
\big| \, $,  $ \, s = \big| \Sigma' \cup \Sigma'' \big| - \big| \Sigma'
\big| - 1 \, $  et  $ \, k = d - \big| \Sigma' \cup \Sigma'' \big| \, $; 
on a, toujours pour les  m\^emes raisons,  $ \, j - \big| \Sigma' \big| <
n - \big| \Sigma' \big| \, $   et  $ \big| \Sigma' \cup \Sigma'' \big| -
\big| \Sigma' \big| - 1 \geq 0 \, $).

\vskip3pt

   En conclusion, on a toujours  $ \, {\big( E' 
\big)}_{\Sigma',\Sigma''} = 0 $,  
d'o\`u  $ \, (E) = 0 \, $,  ce qui termine la preuve.   $ \square $
\enddemo

\vskip23pt


\vskip21pt


\Refs
  \widestnumber\key {WX}

\vskip4pt

\ref
  \key  CP   \by  V. Chari, A. Pressley
  \book  A guide to Quantum Groups
  \publ  Cambridge University Press   \publaddr  Cambridge   \yr  1994
\endref

\vskip2pt

\ref
  \key  Di   \by  J. Dixmier 
  \paper  Introduction to the theory of formal groups   
  \jour  Pure and Applied Mathematics   \vol  20   \yr  1973  
\endref

\vskip2pt

\ref
  \key  Dr   \by  V. G. Drinfel'd
  \paper  Quantum groups
  \inbook  Proc. Intern. Congress of Math. (Berkeley, 1986)  \yr  1987 
  \pages  798--820
\endref

\vskip2pt

\ref
  \key  EK   \by  P. Etingof, D. Kazhdan
  \paper  Quantization of Lie bialgebras, I 
  \jour  Selecta Math. (New Series)   \vol  2   \yr  1996  
  \pages  1--41
\endref

\vskip2pt

\ref
  \key  G1   \by  F. Gavarini 
  \paper  Geometrical Meaning of R--matrix action for Quantum groups 
at Roots of 1 
  \jour  Commun. Math. Phys.   \vol  184   \yr  1997   \pages  95--117
\endref

\vskip2pt

\ref
  \key  G2   \by  F. Gavarini 
  \paper  The  $ R $--matrix  action of untwisted affine
quantum groups at roots of 1 
  \toappear \ in Jour. Pure Appl. Algebra  
\endref

\vskip2pt

\ref
  \key  G3   \by  F. Gavarini 
  \paper  The quantum duality principle  
  \jour  Preprint  
\endref

\vskip2pt

\ref
  \key  On   \by  A. L. Onishchik (Ed.)  
  \paper  Lie Groups and Lie Algebras I   
  \jour  Encyclopaedia of Mathematical Sciences   \vol  20   \yr  1993 
\endref

\vskip2pt

\ref
  \key  Re   \by  N. Reshetikhin 
  \paper  Quasitriangularity of quantum groups at
roots of 1 
  \jour  Commun. Math. Phys.   \vol  170   \yr  1995   \pages  79--99
\endref

\vskip2pt

\ref
  \key  WX   \by  A. Weinstein, P. Xu  
  \paper  Classical Solutions of the Quantum Yang-Baxter Equation  
  \jour  Commun. Math. Phys.   \vol  148   \yr  1992  
  \pages  309--343
\endref

\endRefs
\enddocument